\documentclass{IEEEtran4PSCC}

\usepackage{amssymb}
\usepackage[cmex10]{amsmath}
\usepackage{graphicx}
\graphicspath{{./pic/}}
\usepackage{pstricks}
\usepackage{leftidx}
\usepackage{csquotes}

\newcommand{\nbus}{n_\text{B}}
\newcommand{\nbust}{\widetilde{n_\text{B}}}
\newcommand{\nline}{n_\text{L}}

\newcommand{\Powd}[1]{P_{\text{D},#1}}
\newcommand{\Powg}[1]{P_{\text{G},#1}}
\newcommand{\Pownet}[1]{P_{\text{net},#1}}
\newcommand{\Powline}[1]{P_{\text{L},#1}}
\newcommand{\Powlineall}{P_{\text{L}}}
\newcommand{\Bgl}{B_{\text{G2L}}}
\newcommand{\R}[1]{\mathbb{R}^{#1}}

\newcommand{\PG}{\mathcal{P}_\text{G}}
\newcommand{\PL}{\mathcal{P}_\text{L}}
\newcommand{\PNTC}{\mathcal{P}_\text{NTC}}
\newcommand{\PGt}{\widetilde{\mathcal{P}}_\text{G}}
\newcommand{\PLt}{\widetilde{\mathcal{P}}_\text{L}}
\newcommand{\PNTCt}{\widetilde{\mathcal{P}}_\text{NTC}}
\newcommand{\PFt}{\widetilde{\mathcal{P}}_\text{F}}
\newcommand{\PFat}{\widetilde{\mathcal{P}}_\text{F0}}
\newcommand{\PFatc}{\leftidx{}{\widetilde{\mathcal{P}}}{_\text{F0}^{\mathrm{c}}}}

\newcommand{\PFatic}{\leftidx{}{\widetilde{\mathcal{P}}}{_{\text{F0},i}^{\mathrm{c}}}}
\newcommand{\Abal}{A_\text{balance}}
\newcommand{\Agen}{A_\text{G}}
\newcommand{\bgen}{b_\text{G}}
\newcommand{\Aline}{A_\text{L}}
\newcommand{\Alinei}{A_{\text{L},i}}
\newcommand{\bline}{b_\text{L}}
\newcommand{\blinei}{b_{\text{L},i}}
\newcommand{\Antc}{A_\text{NTC}}
\newcommand{\bntc}{b_\text{NTC}}
\newcommand{\nntc}{n_\text{NTC}}
\newcommand{\Powlinetotal}{P_\text{L}^\text{total}}
\newcommand{\Tline}{T_\text{L}}
\newcommand{\Tbus}{T_\text{B}}

\newcommand{\MapPT}[2]{\mathcal{L}(#1,#2)}

\newcommand{\eref}[1]{(\ref{eq:#1})}
\newcommand{\sref}[1]{Section \ref{sec:#1}}
\newcommand{\fref}[1]{Fig. \ref{fig:#1}}

\begin{document}

\title{Polyhedral Computation Based Transfer Capacities in Multi-Area Power Systems}

\author{
\IEEEauthorblockN{Alexander Fuchs, \emph{Member, IEEE}}
\IEEEauthorblockA{Research Center For Energy Networks \\ ETH Z\"urich, Switzerland\\fuchs@fen.ethz.ch}
\and
\IEEEauthorblockN{Marc Scherer, \emph{Member, IEEE}}
\IEEEauthorblockA{Swiss Transmission System Operator\\Swissgrid Ltd., Switzerland\\marc.scherer@swissgrid.ch}
\and
\IEEEauthorblockN{G\"oran Andersson, \emph{Fellow, IEEE}}
\IEEEauthorblockA{Power System Laboratory\\ ETH Z\"urich, Switzerland\\andersson@eeh.ee.ethz.ch}
}
\maketitle

\begin{abstract}
This paper proposes a new methodology to maximize the feasible set of power injections 
and cross-border power transfers  in meshed multi-area power systems.
The approach used polyhedral computation schemes and is an extension to the classic procedure for cross-border transfer capacity assessment in the European power network, including the computation of bilateral cross-border transfer capacities as well as multilateral flow-based approaches.
The focus is the characterization of inter-area exchange limits 
required for secure power system operation  in the presence of physical transmission constraints, while maximizing the utilization factors of the transmission lines. The numerical examples include a case study of the \mbox{ENTSO-E} transmission system.
\end{abstract}



\section{Introduction}

This paper outlines a methodology to maximize the feasible set of the cross-border transfer capacity calculation. We presents an extension to the classic procedure for cross-border capacity assessment in the European power network. The results are compared to the well-established bilateral computation of cross-border transfer capacities and can be used as an preprocessing step of existing multilateral flow-based approaches.

For the determination of cross-border transfer capacities, the key challenge is the effect of multiple parallel flows through neighboring areas. In the European network, for example, a simultaneous power exchange from Belgium to Italy and from the Netherlands to Austria will induce major flows through the German transmission grid. The transfer capacity must be sufficiently conservative to guarantee the secure system operation during all flow conditions, both regarding the physical constraints of cross-border transmission lines as well as internal transmission constraints of the full networks in the area represented by the nodes.

Our contribution is twofold. First, the capacity bounds consider not only neighbors but all areas between which power exchange is possible, similar to the flow-based methods. Second, the capacity bounds are not a set of box constraints as induced by the upper and lower bilateral transfer capacity limits, but rather a polyhedral capacity set that incorporates the effect of parallel load flows. In the mentioned example, Belgium could likely export much more power to Italy, if there is no export occurring from the Netherlands to Italy, freeing capacities in the German transmission system. Subsequently, the often applied bilateral computation is one specific subset of the proposed polyhedral capacity set.

The paper focuses on the characterization of inter-area exchange limits caused by the physical transmission constraints required for secure power system operation while maximizing the utilization factors of the transmission lines. We investigate neither the allocation of these capacities to the bidding zones nor their pricing.

The structure of the paper is as follows: Section~\ref{European-multi-area-congestion-management} briefly reviews the current European concepts for managing cross-border congestions. Section~\ref{Modeling-of-transfer-capacities} and Section~\ref{Extended-capacity-selection-problem} outline the modeling of transfer capacities and the extended capacity selection problem, respectively. Section~\ref{Simulations} presents the simulation results and discusses the policy implications. Finally, Section~\ref{Closing-remarks} is devoted to conclusions and perspectives.


\section{European Multi-Area Congestion Management}\label{European-multi-area-congestion-management}
In liberalized electricity markets, the transmission system determines the limitation for wholesale and ancillary service trading~\cite{Scherer2013292}. Consequently, the way in which cross-border transfer capacities are calculated has a substantial impact on the market opportunities. In Europe, the Transmission System Operators~(TSOs) bear responsibility for the operation and security of the power system, which includes the determination of the transfer capacity available to market participants' electricity trading. Cross-border capacities are either agreed bilaterally between neighboring countries (\enquote{contract-based}) or determined multilaterally for several areas (\enquote{flow-based})~\cite{OpHB-Team2014,6607366}. Both approaches can be used in regional markets implementations. For example, in Central Western Europe~(CWE), i.e. the Netherlands, Belgium, France, Luxembourg, and Germany, a co-ordinated contract-based market coupling was launched by \mbox{end-2010}, which was superseded by a flow-based day-ahead market coupling \mbox{mid-2015}. We briefly outline the basics of both state-of-the-art approaches, as our methodology to determine the feasible sets for the transfer capacities can be compared to existing bilaterally agreed values as well as enable flow-based transfer capacity calculations.

\subsection{The Available Transfer Capacity Approach}
This mechanism is based on bilateral agreements between neighboring TSOs. Based on historical data, i.e. reference days, well-representing seasonal patterns as well as justified security margins, each TSO determines a Total Transfer Capacity~(TTC) for each direction on each border of its control area. Thus, the TTC is the upper limit for which the maximum physical flow on a critical network element does not exceed its safety margins, i.e. \mbox{N-1} criterion. Based on the TTC, each TSO deducts a safety margin, referred to as the Transmission Reliability Margin~(TRM), as well as used capacities of Long-Term Contracts~(LTCs), as a holder of a LTC must always declare by the previous day, whether or not and to what extent the holder intends to use the respective long-term reserved transfer capacities. The Net Transfer Capacity (NTC) available to wholesale trading results from the TTC minus the TRM minus the LTC. Finally, the Available Transmission Capacity~(ATC) is the part of NTC that remains available after each phase of the allocation procedure for further commercial activity, i.e. ATC is NTC minus Already Allocated Capacity~(AAC). This whole process is operated bilaterally; therefore, if the calculated values deviate between neighboring TSOs, generally the lower ones are selected. By that, borders are considered separately which does not allow for a holistic consideration of the power flows in the power system.

\subsection{The Flow-Based Approach}
Instead of fixed capacities, the flow-based methodology is based on a reduced network constituting of nodes and lines in order to take into account that electricity can flow via different paths in an highly meshed power system. Each TSO provides input data, which is combined at a regional level. To retrieve a reduced network model, instead of considering each and every line, so-called critical branches are introduced. They consist of tie-lines as well as internal lines that significantly impact a cross-border exchange. This allows to determine which combinations of cross-zonal exchanges may lead to an overload of a critical network element. Based on the physical limit and potential security margins of the line, the physical capacity for each critical branch, i.e. the total maximal flow, is determined. A Flow Reliability Margin~(FRM) is to cope with the uncertainty inherent to the process of determining the remaining capacity, and a reference flow to consider the already known long-term nominations. What will eventually be offered to the wholesale market is the so-called Remaining Available Margin~(RAM). By that, it is possible to consider the effect of a meshed system.

\section{Modeling of transfer capacities}\label{Modeling-of-transfer-capacities}

This section introduces the mathematical modeling needed for the characterization of transfer capacities.
The modeling is the same, whether it is applied to characterize TTC, NTC, ATC or the flow-based approach.
The first part defines the full power system model and constraint formulation.
The second part defines the aggregation scheme used to obtain a reduced network model.

\subsection{Power system modeling and constraints}
\label{sec:psmodel}
The power system is modeled using a DC power flow approximation  \cite{MachPow}.
The  network has $\nbus$ buses connected by $\nline$ lines.
Each bus $i$ has a known demand power $\Powd{i}$ and a generator power  $\Powg{i}$ 
that is to be selected during the power flow optimization.
The difference between generation and demand at each bus is the net power injection
\begin{equation}
\Pownet{i}  = \Powg{i} - \Powd{i}, \qquad i = 1,2, ..., \nbus \quad , 
\end{equation}
and forms the power system state 
$x \in \R{\nbus}$, 
    \begin{equation}
x = [\Pownet{1}, \Pownet{2}, ...,  \Pownet{\nbus} ]^T \quad .
    \end{equation}
Furthermore, the vector of active power flows in the $\nline$ lines is defined as $\Powlineall \in \R{\nline}$, 
    \begin{equation}
\Powlineall = [\Powline{1}, \Powline{2}, ...,  \Powline{\nline} ]^T \quad .
    \end{equation}
 In the DC power flow model, all relations between bus voltage angles, net power injections and power flows in the lines become linear, in particular \cite{matpower}
\begin{equation}
\Powlineall =  \Bgl x \quad .
\label{eq:Bgl}
\end{equation}
The net power injections have to satisfy the power balance equality constraint, 
\begin{equation}
 \Abal^T x = [1, 1, ..., 1] x =  \sum_{i=1}^{\nbus} \Pownet{i}  = 0\quad.
\end{equation}
There are two sets of inequalities constraining the set of possible net power injections.
The generator powers at each node are bounded by the generation limit
\begin{equation}
0 \leq \Powg{i} \leq  \overline{\Powg{i}} \quad , 
  \label{eq:pgcon}
\end{equation}
and the power flows in each  line are bounded by  the thermal limit
\begin{equation}
 |\Powline{i}| \leq  \overline{\Powline{i}} \quad .
  \label{eq:plcon}
\end{equation}
The inequalities define the generator constraint polyhedron $\PG$ and the line constraint polyhedron $\PL$,  
\begin{align}
\PG & = \{x \in \R{\nbus} \ : \quad  \Agen x \leq \bgen, \ \Abal x = 0 \} \label{eq:PGpoly}\\
\PL & = \{x \in \R{\nbus} \ : \quad  \Aline x \leq \bline, \ \Abal x = 0 \}, 
\end{align}
where the parameters $\{\Agen, \bgen\}$  are computed from \eref{pgcon} and
the parameters $\{\Aline, \bline\}$ are computed from \eref{Bgl} and \eref{plcon}.

A system state $x$ that satisfies all inequalities lies in both polyhedra,
  \begin{equation}
x \in (\PG \cap \PL)\quad , 
  \end{equation}
and is referred to as \emph{feasible}, otherwise it is referred to as \emph{infeasible}.

\subsection{Network aggregation}

The background of this paper is the aggregation of a detailed power system model as defined in \sref{psmodel} with the system state $x \in \R{\nbus}$ 
by assigning the original $\nbus$ buses to one of the $\nbust$  regions, with $\nbust < \nbus$.
The outcome is a reduced power system model with the system state $y \in \R{\nbust}$ denoting the total net power injections of the regions.
The transformation with the  \emph{bus aggregation matrix} $\Tbus$
\begin{equation}
y = \Tbus  x
\label{eq:yTx}
\end{equation}
is essentially a summation of the net power injections of all buses associated with a region.
If the original bus $i$ is assigned to region $j$, then the element of $\Tbus$ in column $i$ and row $j$ is 1, and zero otherwise.

The linear map \eref{yTx} also defines a mapping $\MapPT{\cdot}{\cdot}$ of polyhedral sets from the original to the reduced state space.
For instance, 
    \begin{align}
\PGt & = \MapPT{\PG}{\Tbus} \nonumber \\ 
     & = \{y\in\R{\nbust} :  \exists x \in \PG: y = \Tbus x \}
    \end{align}
denotes all reduced system states with a corresponding original state that satisfies the generator constraints.
Since $\nbust < \nbus$, $\PGt$ is a projection on the dimensions defined by the rows of the bus aggregation matrix $\Tbus$.
Linear mappings of polyhedra can be computed exactly or with approximations using existing software implementations \cite{mpt3}, \cite{colinproj}.

The aggregation transformation is assumed to be known and fixed. 
It originates from the organization of the power system into control areas or market zones.

\subsection{Feasibility conditions}

To characterize feasible power injections in the reduced and original system,  two conditions are defined.

Given a generation constraint polyhedron $\PG$, a line constraints polyhedron $\PL$ and a bus aggregation matrix $\Tbus$, 
a state $y$ of the reduced system is \emph{feasible in the  original system} if 
\begin{equation}
\exists x \in\PG :   y = \Tbus x, \quad x \in \PL \quad .
\label{eq:feascon}
\end{equation}
For a reduced system state $y$ that satisfies \eref{feascon}, it is by definition always possible to find a corresponding power injection $x$ that satisfies all generator and line constraints. 
However, the identification of the corresponding power injections is not unique. 
In fact, some of the power injections that satisfy the generator constraints may violate internal line constraints and have to be avoided through a regional dispatch.

Consequently,  $y$ is referred to as \emph{strongly feasible in the  original system} if it satisfies \eref{feascon} and additionally  
\begin{equation}
\nexists x \in\PG :   y = \Tbus x,  x \notin\PL \quad. 
\label{eq:strongfeascon}
\end{equation}
In this case, a regional dispatch is not required since the line constraints are automatically satisfied with the generator constraints.

\section{Characterization of feasible network injections}\label{Extended-capacity-selection-problem}

This section characterizes the original power system constraints in the reduced power system model.
      It enables operational decisions with the reduced model that automatically ensure the absence of constraint violations in the original system.
The required information for decisions is only the reduced system state $y$,  the original system model  and state are not needed.

Three polyhedral constraint sets of the reduced system state $y$ are presented,  that can be calculated in a preprocessing step from the original system model.

The first two sets, parameterizing the NTC approach and the flow-based approach,  ensure the feasibility condition \eref{feascon}. 
The third set ensures the strong feasibility condition \eref{strongfeascon}.
All sets are optimal in the sense that they maximize the transmission line utilization possible with the approach.

\subsection{Mapping of NTC constraints }
\label{sec:ntcmap}
NTC bounds constrain the aggregated active power flow $\Powlinetotal$ through selected transmission lines of the network,  
\begin{equation}
\Powlinetotal= \Tline \Powlineall \leq \bntc\quad .
\label{eq:Pntctotal}
\end{equation}
If the original bus $\Powline{i}$ contributes to the $j$'th NTC bound of the NTC vector $\bntc$, then the element of the line aggregation matrix $\Tline$ in column $i$ and row $j$ is 1, 
   and zero otherwise.
         
The polyhedron of power injections satisfying the NTC and generation constraints is defined as 
\begin{align}
\nonumber
\PNTC = \{x\in\R{\nbus}:  \Antc x \leq \bntc ,   \Agen x & \leq \bgen, \\
         \Abal x &  = 0 \ \ \}  \quad,
\end{align}
with $\Antc = \Tline \Bgl $ combining \eref{Bgl} and \eref{Pntctotal} and the other parameter as in \eref{PGpoly}.
The mapped NTC polyhedron
\begin{equation}
\PNTCt = \MapPT{\PNTC}{\Tbus}
\end{equation}
characterizes the NTC constraints in the reduced system space.

The states in the mapped NTC polyhedron can satisfy the feasibility condition \eref{feascon} if the NTC bounds $\bntc$ is not chosen too large.
The maximum NTC bounds that still satisfy the feasibility condition are obtained by solving
\begin{equation}
\max_{\bntc}\quad w^T \bntc \quad  \text{s.t.} \quad  \PNTC \subset (\PG \cap \PL)
\label{eq:ntcopt}
\end{equation}
with the vector $w$ denoting a weighting of the NTC bounds.
Since the polyhedral set constraint in \eref{ntcopt} complicates the problem, 
      it is useful to simplify the problem by searching along a nominal NTC directions $\nntc$, 
\begin{align}
 k_i^* = \min_{k_i, x}\quad  & k_i \label{eq:linntc1}\\ 
 \text{s.t.} \quad  \Antc x &  \leq \nntc\cdot k_i ,\\
                                            \Agen x &  \leq \bgen, \\
                     \label{eq:linntcend}                       \Alinei^T x & \geq \blinei, 
\end{align}
where $\Alinei^T$ denotes the $i$'th row of the line constraint matrix $\Aline$. 
The solution of \eref{linntc1}-\eref{linntcend} determines the smallest NTC scaling along direction $\nntc$ that violates the $i$'th line constraint.
Repeating the linear program for all line constraints yields the optimal NTC scaling
\begin{equation}
k^* = \min_{i }\quad k_i^* \quad .
\end{equation}
with $\bntc = k^*\nntc$. 
The procedure can be repeated for randomly sampled NTC directions to maximize a common objective like the sum of all NTCs or an economic objective.
If an NTC direction and scaling is provided as parameter, for instance from bilateral agreements,  the linear programs can still be solved without an objective as a pure feasibility problem to verify that the NTC bounds prevent all constraint violations.

During operation, reduced system states in  $\PNTCt$ with $\bntc$ selected using the linear program, satisfy  a sufficient condition for the feasibility condition \eref{feascon}.
The evaluation of a candidate system state $y$ is very fast, requiring only the computation of $\Powlinetotal$ through a matrix multiplication of $y$ and 
the verification of the box constraints \eref{Pntctotal}.

\subsection{Mapping of full flow constraints}
\label{sec:flowmap}

The necessary and sufficient for the feasibility condition \eref{feascon} define a set containing \emph{exactly all} reduced states $y$ that have a corresponding feasible original state $x$.
The computation of this set is a key preprocessing step of the flow-based dispatch since it allows the maximum utilization of the available transmission capacities.
Mathematically, the set is a linear map of the feasible constraint polyhedron to the reduced system state, 
\begin{equation}
\PLt= \MapPT{\PG \cap \PL}{\Tbus} \quad .
\end{equation}

The complexity of the projection mainly depends on the dimension of the reduced space and the number of line constraints. 
Inner approximations of the projected set can be efficiently obtained using sampling approaches similar to the NTC maximization, since the constraint sets are bounded intersections of convex sets. 

During operation, the evaluation of a candidate system state $y$ is  still very fast, requiring only the evaluation of the polyhedral  set constraints of $\PLt$ through a matrix multiplication and vector comparison.

\subsection{Set difference of line violations }
The reduced system states $y$ satisfying the strong feasibility conditions form a subset of $\PLt$. 
In addition to \eref{feascon}, they have to satisfy the constraint \eref{strongfeascon},   which defines a set
\begin{align}
\PFat & = \{ y \in \R{\nbust}: 
\nexists x \in\PG :   y = \Tbus x,  x \notin\PL 
\} \quad ,\\
 & = \{ y \in \R{\nbust}: 
\exists x \in\PG :   y = \Tbus x,  x \notin\PL 
\}{^{\mathrm{c}}} \quad,  \\
        & = \left(\PFatc\right){^{\mathrm{c}}}
\end{align}
written in the second line using the complement.
The set $\PFat^c$ is a collection of overlapping polyhedra, with 
\begin{align}
 \PFatc  & = \cup_i \  \PFatic\\
&  =  \cup_i \ 
  \{  y \in \R{\nbust}: 
\exists x \in\PG :   &&y   = \Tbus x, \nonumber \\
                     & &&      \Alinei^Tx  \geq\blinei 
 \} \quad  . \quad  \label{eq:PFatviolate}
\end{align}
The line constraint violation in \eref{PFatviolate} covers all cases where 
$x \notin \PL$.
The final step is the computation of the set $\PFt$ that satisfies both conditions \eref{feascon}-\eref{strongfeascon}, resulting in 
    \begin{align}
\PFt & = \PGt  \cap  \left( \PFatc\right){^{\mathrm{c}}}\\
         & = \PGt \backslash   \PFatc \quad.
    \label{eq:PbP}
    \end{align}
The set difference operation in \eref{PbP} of a polyhedron $ \PGt$ with a family of overlapping polyhedra $ \PFatc $ is described in \cite{Mato}.
It successively constructs a family of non-overlapping polytopes covering the non-convex set $\PFt$ and is available as software implementation \cite{mpt3}.

The outlined approach provides a set of  constraints that are necessary and sufficient for  the strong feasibility condition.
The resulting set $\PFt$ enables the selection of reduced system states $y$ that are guaranteed to satisfy the line constraints of the original model, no matter what generator settings are selected.
This is particularly useful if during the power system dispatch the actual 
distribution of generator capacity in each region is not known, for example due to fluctuating availability of renewable energy sources or restricted communication.

On the downside, the approach is more complex than the other two approaches.
For the preprocessing,  a polyhedral set difference computation is required.
The evaluation of a candidate system state $y$, requires the evaluation of potentially many polyhedra defining the non-convex region.
Several approaches to reduce the evaluation complexity exist, but require addition preprocessing \cite{Tondel}, \cite{optitrees}, \cite{LDF}.

Finally, the strong feasibility condition can quickly become too strict, since some infeasible power system states map to the same reduced system state,  resulting in $\PFt$ being the empty set. 
A solution can be a more detailed modeling, including the violated transmission lines in the reduced system model.
Alternatively, only the regular feasibility condition can be required as in \sref{ntcmap} and \sref{flowmap}, thereby shifting the responsibility for the local line constraints to the regional grid participants. 

\section{Numerical examples}\label{Simulations}
This section applies the proposed methodology to two example systems. 
A simple system highlights the different feasibility sets. 
Then, a large model of the ENTSO-E network is used to demonstrate the applicability of the method to practical system models.

\subsection{Illustrative example}
\label{sec:7bus}

The simple example system consists of six buses that are aggregated into three different areas.
The topology, line constraints  and power capacity is given in \fref{sixbus}.
        All lines are 500 km long, using an inductance of $0.09 \text{ p.u.}/\text{km}$ for a base voltage of 380 kV and a base power of 900 MVA. 
        The southern part of the system is a net power importer from the center or the north. 
        The limiting factor is caused by the constrained transmission inside the central region.
\begin{figure}
    \center
    \includegraphics[width=.5\columnwidth]{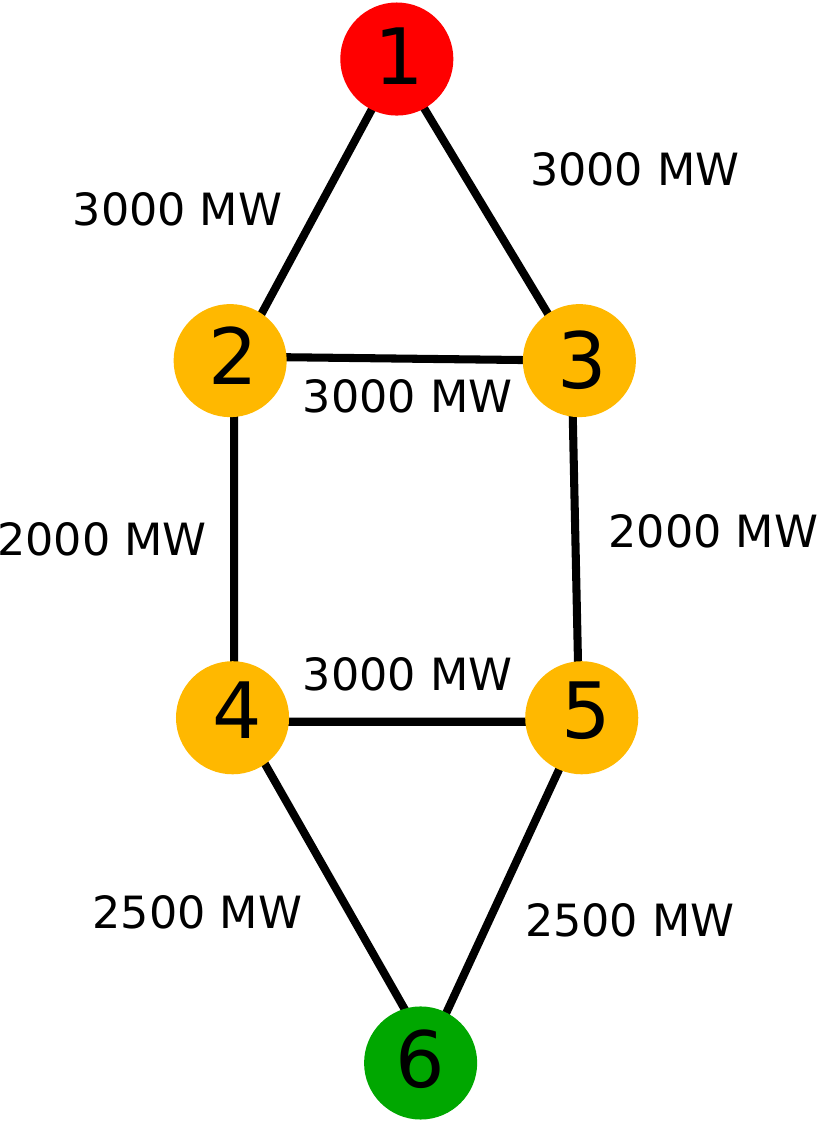}
    \caption{Topology of the six bus example system, showing bus numbers and line constraints.
        The nodes 2-5 are aggregated into the center region.
   The demand power at node 6 is 6 GW, at all other nodes 3 GW. 
   The maximum generation power at nodes 1 and 4 is 10 GW, at all other nodes 3 GW.
   1-5 have a demand of 3 GW, node 6 has a demand of 6 GW.
    }
    \label{fig:sixbus}
\end{figure}

After the aggregation step, only three net power injections remain and $y\in\R{3}$.
Furthermore, one of these free variables becomes dependent on the others due to the power balance constraint. 
Therefore,  the full feasible set can be illustrated in a 2D plot without any further projections being required.
The resulting polytopes can also be plotted in $\R{3}$ and will lie on a hyperplane in the $[1, 1, 1]^T$ direction.

\fref{7bus_ntc} shows that the full projection of the line constraints yields a much larger feasible set than the largest set obtainable with pure NTC constraints between the three regions. 
The reason is that basically a larger power export in the north is possible whenever the center reduces its power production. This effect requires the flow-based approach that takes into account the overall state of the network, not only isolated lines.

\begin{figure}
    \center
    \def\svgwidth{.9\columnwidth}
    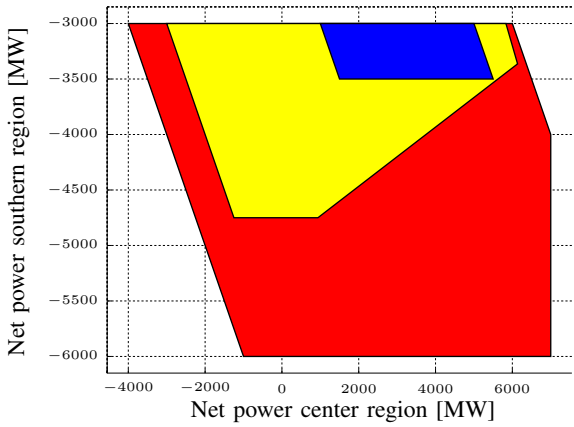
    \caption{Six bus system aggregated into three regions.  Possible net power injections in the center and southern region.
        Generator constraints $\PGt$ (red), 
    line constraints $\PLt$ (yellow) and 
    NTC constraints $\PNTCt$ (blue).
    }
    \label{fig:7bus_ntc}
\end{figure}

\fref{7bus_strong} illustrates the strong feasibility approach. 
The goal is to determine the subset of $\PLt$ that not only has \emph{some} feasible corresponding state in the original system, but for which \emph{all} the corresponding states that can be realized by the generators do not violate any line constraints.
This requirement is quite strict an required to relax the line constraints of the six bus system. 
To this end, the capacity of the two transmission lines in the center are increased to 2700 MW and the two souther line capacities to 3300 MW. 
The result shows the non-convexity of the strongly feasible set, caused by the set difference operation.

\begin{figure}
    \center
    \def\svgwidth{.9\columnwidth}
    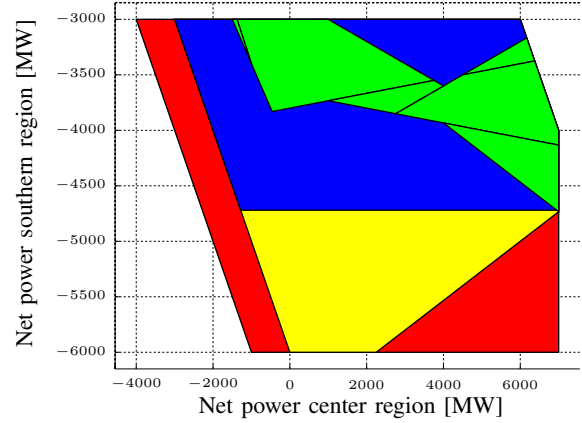
    \caption{Six bus system with increased transmission capacity aggregated into three regions. Possible net power injections in the center and southern region.
        Generator constraints $\PGt$ (red), 
    line constraints $\PLt$ (yellow), 
    NTC constraints $\PNTCt$ (blue) and 
    strong feasibility constraints $\PFt$ (green). 
    }
    \label{fig:7bus_strong}
\end{figure}

\subsection{European system}
\label{sec:entso}
The method is now applied to a large model of the ENTSO-E transmission grid, documented in \cite{pegase} with an implementation available in \cite{matpower}.
It consists of 9241 buses, 16049 lines and is aggregated into 23 zones as illustrated in  \fref{entsoemap}.
The resulting system can now be analysed regarding the constraints of the individual zones or larger subgroups of the system.

Two results are shown in this section.
\fref{polish_tope} illustrates the constraint coupling between the German and the Polish zone.
Note that the figure shows only a projection of the feasibility set onto these two dimensions.
A specific flow pattern between the other control zones is required to realize the maximum possible export of 10 GW from the polish zone. 
Using the classical NTC or ATC characterization as outlined in Section II
and \cite{ntcdef}, simplifies the constraint set but makes it also more conservative. 
The illustration applies the Polish NTC values from 2011 \cite{ntcmatrix} to the different neighboring lines with a total of 3500 MW of cross-border transfer capacity.
The violation observed along the German dimension is no issue, since additional NTC constraints of other critical system boundaries are commonly used to ensure full feasibility.

     Finally, \fref{entsoe_tope} shows the application of the proposed method 
     to determine the full constraint set if the ENTSO-E network is operated in four distinct areas.
     The grouping into the four regions is also shown in \fref{entsoemap}.
    For the constraint characterisation, all of the 16049 original line constraint are included. 
    The projection step uses a sampling step to determine the  characteristic directions in the reduced network space.
The accuracy is ensured by comparing inner and outer approximations of the constraint set.     

   The projection result is a polytope in the four dimensions corresponding to the four network regions.
   Due to the intersection with a hyperplane representing the power balance constraint, it is possible to display the constraint set in a 3D picture.
For the application, also groupings into more than 4 regions can be considered and add little to the evaluation complexity.

\begin{figure}
    \center
    \includegraphics[width=.9\columnwidth]{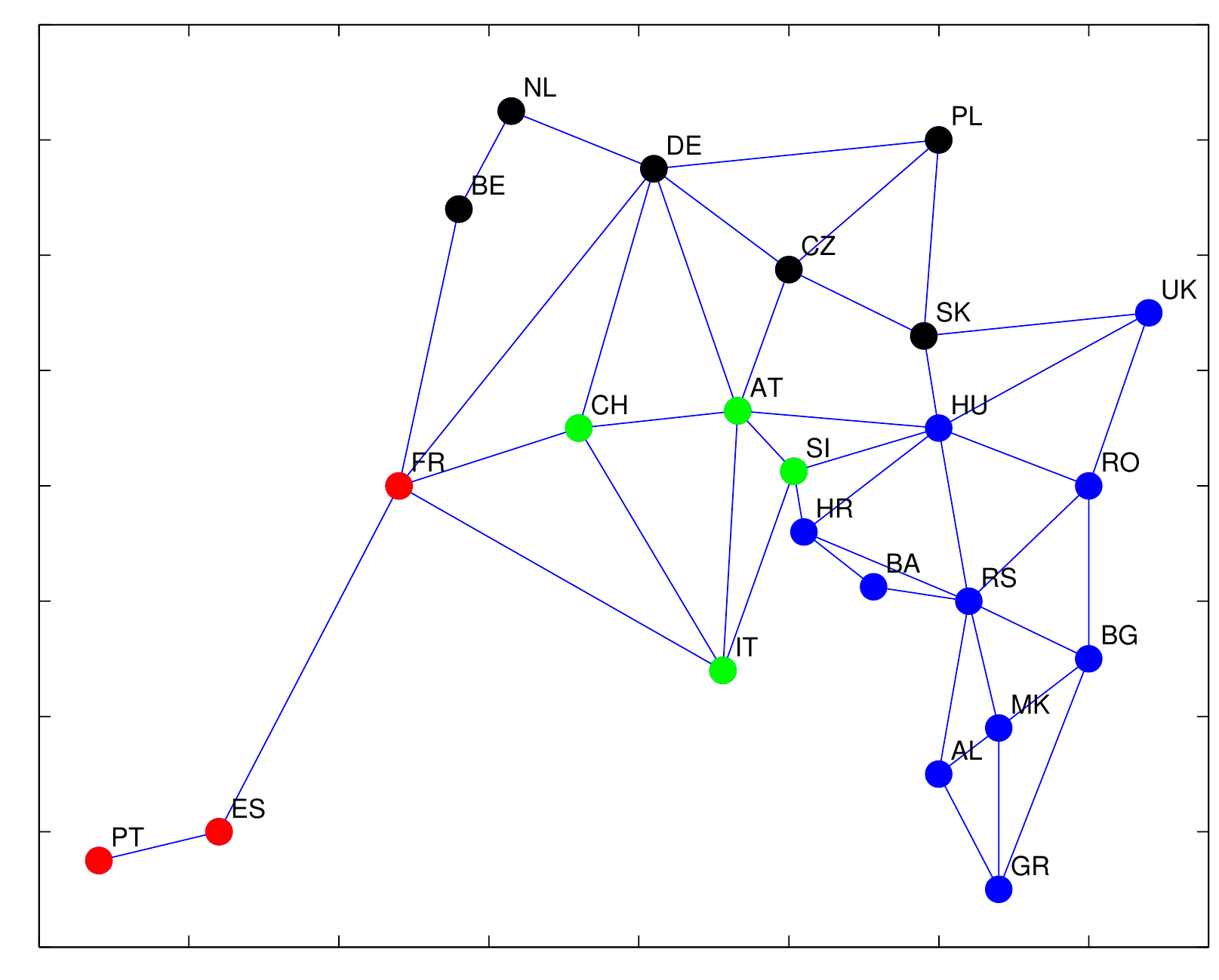}
    \caption{Map of the ENTSO-E system with 9241 buses, aggregated in to 23 zones (circles) according to \cite{pegase}, 
        and into four large regions to illustrate cross border power exchanges: 
    South-western (red), northern (black),  central (green) and south-eastern region (blue).
    }
    \label{fig:entsoemap}
\end{figure}

\begin{figure}
    \center
    \def\svgwidth{\columnwidth}
    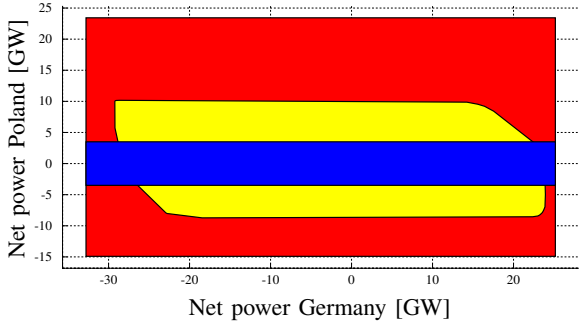
    \caption{Possible net power injections of the Polish zone and the German zone, projected from the full ENTSO-E system model.
Feasible set defined by generator constraints $\PGt$ (red), 
  the full  line constraints $\PLt$ (yellow) and the Polish NTC constraints $\PNTCt$ (blue) according to \cite{ntcmatrix}.
       Note that relying purely on classical NTC definitions significantly underestimates the admissible region of net power injections.
    }
    \label{fig:polish_tope}
\end{figure}

\begin{figure}
    \center
    \def\svgwidth{\columnwidth}
    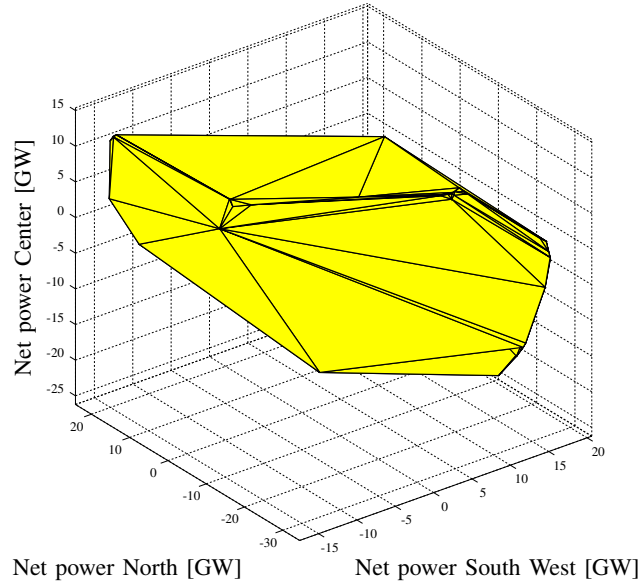
    \caption{Possible net power injections of the ENTSO-E system aggregated into four regions, projection of the line constraints $\PLt$.
         The feasible set of the south-western, northern and central region is shown.
        The south-eastern region is not shown and compensates the total power imbalance of the shown regions to zero.
        The polyhedral representations of admissible net power injections allows a fast feasibility assessment of the aggregated system model.
    }
    \label{fig:entsoe_tope}
\end{figure}


\section{Closing Remarks}\label{Closing-remarks}
\subsection{Conclusion}
This  paper demonstrated how the  constraints of a detailed large-scale network model can be characterized in a reduced aggregated power system model.
The proposed method enables the operator to take decisions in the reduced model, that automatically ensure the absence of constraint violations in the detailed network model. 
We outlined different polyhedral constraint sets that can be calculated in a preprocessing step of the market-driven transfer capacity allocation. 
The method was evaluated using a six bus test system as well as a large-scale network model of \mbox{ENTSO-E}. 
The results show  a large potential of the flow-based capacity allocation if the full constraint set can be used, compared to a classical NTC approach.

\subsection{Outlook}
The next key step is the application of this method to  a market model with a realistic bidding process.
The European Union's "Third Energy Package" stipulates a competitive and integrated European electricity market with extensive cross-border trade facilitation.
 This not only emphasizes the need for transparent and traceable methods to determine cross-border capacities, it also requires system operators and respective regulators to ensure the effective and optimal determination of those. For a large-scale implementation, the practical applicability of the proposed method needs to be further investigated.
In this paper, the reduced  network model is assumed to be given from organizational or political constraints. 
The choice of the appropriate aggregation scheme and reduced model has a large impact on the constraint sets as well as computational efforts during the decision making 
and is therefore another topic of investigation.


\bibliographystyle{IEEEtran}
\bibliography{fuchs-scherer_PSCC2016}

\end{document}